\documentclass[twocolumn,7pt]{article} %% use 11pt for Applied Optics
\usepackage{amsmath,amssymb,graphicx}
\usepackage[hyperindex,breaklinks]{hyperref}
\hypersetup{colorlinks, citecolor=blue, linkcolor=blue, urlcolor=blue}
\usepackage{csquotes}
\usepackage{color}

\newcommand{\dd}{\;\text{d}}

\usepackage{graphicx}
\usepackage[a4paper]{geometry}
\usepackage[affil-it]{authblk}
\usepackage{scrextend}

\makeatletter
\newcommand\footnoteref[1]{\protected@xdef\@thefnmark{\ref{#1}}\@footnotemark}
\makeatother

\newcommand{\e}{\text{e}}

\usepackage{placeins}

\begin{document}

\title{\textbf{Noise Robustness of a Combined Phase Retrieval and Reconstruction Method for Phase-Contrast Tomography}}

\author{ Rasmus Dalgas Rasmussen\footnote{Department of Applied Mathematics and Computer Science, Technical University of Denmark.} ,  Jakob Sauer J{\o}rgensen$^{*}$, Henning Friis Poulsen\footnote{Department of Physics, Technical University of Denmark.} \hspace{0.7mm} and \hspace*{1.7cm}  Per Christian Hansen$^*$} %\affil{DTU Compute \and DTU Compute}
%\author{Jakob Sauer J{\o}rgensen} \affil{DTU Compute}

\maketitle

\begin{abstract}
Classical reconstruction methods for phase-contrast tomography consist of two stages: phase retrieval and tomographic reconstruction. 
A novel algebraic method combining the two was suggested by Kostenko et al.\ (Opt.\ Express, 21, 12185, 2013) and preliminary results demonstrating improved reconstruction com-
pared to a two-stage method given. 
Using simulated free-space propagation experiments with a single
sample-detector distance, 
we thoroughly compare the novel method with the two-stage method to address limitations of the preliminary results. 
We demonstrate that the novel method is 
substantially more robust towards noise; our simulations point to a possible reduction in counting times by an order
of magnitude.
\end{abstract}

\section{Introduction}

With the upcoming of high-brilliance X-ray sources, new phase-contrast tomography (PCT) methods have gained widespread use \cite{Bravin2013}. Among these is the free-space propagation method \cite{Cloetens1996,Nugent1996a}, with the advantage that no additional optical elements such as analyzer crystals or gratings are required. Compared to conventional absorption contrast, phase contrast may provide adequate contrast at lower dose rates, thus allowing segmentation of objects comprising two or more materials with nearly the same electron density (e.g., \cite{Salvo2003}). For a comparison of variants of free-space propagation methods in general see \cite{langer2008}.

First experiments were performed using holo-tomography, requiring the combination of measurements from several sample-to-detector distances \cite{Cloetens1999}, but today highly successful reconstructions are often possible using a single distance, e.g., in the area of paleontology \cite{Tafforeau2006a}. This is experimentally convenient, but also remarkable, as the information content is obviously reduced. In fact, reconstruction with the single-distance set-up is typically based on the work by Paganin et al.~\cite{Paganin2002}, assuming proportionality between the absorption and the phase shift. Research into the limitations of this so-called duality method will benefit experimental planning.

The standard approach to PCT reconstruction is a two-stage procedure. In the first stage the phase and absorption fields are determined for each projection using a phase retrieval algorithm. In the second stage, a classical algorithm is used to compute a reconstruction based on the projection fields.

Recently, Kostenko et al. proposed a combined approach \cite{Kostenko2013a,Kostenko2013} %that could extend the range of applicability of the duality method.
%with the possible advantage of suppressing mathematical modelling errors through regularization.
%We see this as encouraging since such an approach could make the duality method more robust and possibly extend its range of applicability.
%From Kostenko's work the combined method
which -- in several of their simulated experiments with noise-free data -- performs better in terms of reconstruction error.
However, for simulated noisy data the combined method is outperformed by the two-stage method. Both methods are tested on simulated data with artificial material indices,
 %but for real data only the combined method is tested in \cite{Kostenko2013a}.
and therefore it is unclear if the combined method performs better than the two-stage method for realistic samples, for different material types, and for varying noise levels.

In this paper we provide a careful numerical simulation study of the combined duality method, comparing it to the two-stage approach.
The use of regularization is key to obtaining high-quality reconstruction, and we
focus on the use of total variation (TV) \cite{Rudin1992} as the regularizer.
%, since the use regularization is key to obtaining high-quality reconstruction, and we use a total variation regularization method \cite{Rudin1992}.
Many samples in materials science and geoscience comprise discrete objects (grains, fibres, cracks), and edges naturally lead to high phase contrast in the free-space propagation method. Using a polycrystal with small density variations as a phantom and with simulated noise, our simulations aim to carefully compare the reconstruction capabilities of the two methods.
%determine difference in performance with respect to small density variations and noise.

After a review of a classical PCT reconstruction method and the combined method, we present our numerical implementation. %where we use state-of-the-art methods. 
In simulations with realistic material parameters we perform a comprehensive comparison of the two methods with respect to different material parameters of increasing difficulty and to the robustness towards noise.

Our results show that the combined method produces improved reconstructions across the range of low to high-absorption materials with small absorption contrast. Furthermore, as the simulated noise level is increased, reconstructions from the combined method show much greater robustness to the noise. This could be of critical importance in practical applications where noise is always a concern.
\section{Definitions and forward model}

In this section we briefly review the underlying definitions and models.
To simplify the presentation and the numerical experiments we consider only 2D problems.

Scalar functions are denoted with italic, e.g., $u$, $B$ or $\phi$. Vectors are denoted with subscript $v$, e.g., $x_v$ or $B_v$. Matrices are denoted with bold upper case, e.g., $\mathbf{A}$ or $\mathbf{F}$, and operators with upper case calligraphic letters, e.g., $\mathcal{A}$ and $\mathcal{F}$. Throughout, $\|\cdot \|$ denotes the vector 2-norm.

\subsection{Definitions}

The Fourier transform $\mathcal{F}$ of a one-dimensional signal $f$ is denoted with $\widehat{\cdot}$:
  \begin{equation}
    \widehat{f}(\omega) = \mathcal{F}(f(x) ) = \int_{-\infty}^\infty f(x) \e^{-2\pi \hat{\iota} \omega x} \dd x.
  \end{equation}
The independent frequency variable is $\omega$ and the complex unit is denoted as $\hat{\iota}=\sqrt{-1}$. The Radon transform $\mathcal{R}$ of a two-dimensional signal $f(x_1,x_2)$ is defined by:
\begin{align}
&[\mathcal{R}f] (t,\theta) = \nonumber \\
&\int_{-\infty}^\infty f\Big(t\cos(\theta)-\tau\sin(\theta),t\sin(\theta)+\tau\cos(\theta)\Big) \dd \tau .
\end{align}
Here $\theta$ is the angular variable and $t$ is the translational variable, perpendicular to line integration direction, in other words the coordinate variable on the one-dimensional detector.

A discrete linear inverse problem can be formulated for absorption-based computed tomography (CT) by discretizing a 2D object into $N$ by $N$ square pixels with pixel values stacked into a vector $u_v$. The Radon transform is discretized using the line-intersection method and represented using a system matrix $\mathbf{A}$ with elements $a_{mn}$ of the path length of X-ray $m$ through pixel $n$. Letting $b_v$ denote the discrete projection data, the discrete linear inverse problem can be written
\begin{align} \label{eq:DiRadon}
\mathbf{A} u_v = b_v.
\end{align}
For more details see \cite{Buzug2008}.

\subsection{Free-space propagation model}

\begin{figure}[htbp]
\begin{minipage}{0.18\columnwidth}
\centering
{\large Coherent wave field } \\[0.3cm]
{\Huge $\Longrightarrow$}
\vspace{1.6cm}
\end{minipage}
\begin{minipage}{0.78\columnwidth}
\includegraphics[trim = 41.5mm 6mm 7mm 0mm, clip, width=\linewidth]{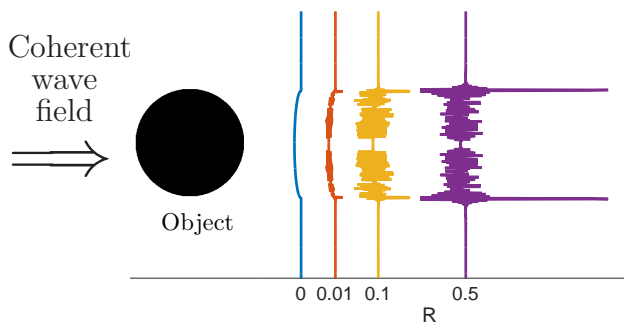}
\end{minipage}
\caption{2D sketch of free-space propagation set-up. The intensity profile is shown as function of increasing detector distance $R$. In the contact plane $R=0$ only absorption contrast is visible. In the near-field region, of interest here, both absorption and phase contrast contribute to the intensity.}
\label{fig:PCT2D3}
\end{figure}

Different experimental PCT setups exist, see \cite{Bravin2013} for an overview. In the present work we focus on the free-space propagation method, \autoref{fig:PCT2D3}, which in some sense is the simplest since it does not require analyser crystals, gratings or likewise. The set-up is similar to the standard CT scanning set-up, with the added requirement that the X-ray source is partly coherent.

%\figone{2D_Sketch}{1}{2D sketch of free-space propagation set-up. The intensity profile is shown as function of increasing detector distance. In the contact plane only absorption contrast is visible. In the near-field region, of interest here, both absorption and phase contrast contribute to the image. Taken from \cite{Bronnikov2002}.}

In free-space propagation PCT phase shifts of X-ray waves are magnified as the object-to-detector distance increases. Based on intensity measurements at the detector, the task is to reconstruct both the absorption index and the refractive index decrement of the object.

For free-space propagation PCT the measured intensity data in the Fresnel region is related to the material properties through a non-linear propagation model. For a 2D object with absorption index $\delta(x_1,x_2)$ and refractive index decrement $\beta(x_1,x_2)$, the corresponding projections $B(t,\theta)$ and $\phi(t,\theta)$, here called \textit{absorption} and \textit{phase shift}, can be modelled as line integrals along the X-ray propagation \cite{langer2008}:
\begin{align}
\label{eq:Btheta}
B_\theta(t)&=\phantom{-}\frac{2\pi}{\lambda} [\mathcal{R}\beta](t,\theta), \\ %= \int \beta(\mathbf{x}) \dd \mathbf{x},
\phi_\theta(t) &=\, -\frac{2\pi}{\lambda} [\mathcal{R}\delta](t,\theta). %\int \delta(\mathbf{x}) \dd \mathbf{x}.
\end{align}
Here $\lambda$ is the X-ray wavelength and $\mathcal{R}$ is the Radon transform. Based on the absorption and phase shift, the measured intensity $I_\theta^R$ is modelled as the squared absolute value of the convolution of the transmittance $T_\theta$ and the Fresnel propagator $P^R$:
\begin{align}
T_\theta(t) &= \exp\Big(-B_\theta(t) + \hat{\iota}\phi_\theta(t)\Big), \label{eqn:transmit} \\
P^R(t) &= -\tfrac{\hat{\iota}}{\lambda R}\exp\Big(\tfrac{\hat{\iota}\pi}{\lambda R}|t|^2 \Big), \\
I^R_\theta(t) &= |T_\theta(t)\star P^R(t)|^2. \label{eqn:intensity}
\end{align}
Here $R$ is the object-to-detector distance, $|\cdot|$ the absolute value and $\star$ the convolution operator.

Discretization of the domain and hence the object of interest into $N \times N$ pixels, gives us discrete versions of the material parameters $\beta$ and $\delta$. A single index notation is introduced $n=i+(N-1)j$, for $i,j = 1,2,...,N$, so $n=1,2,...,N^2$, which gives two column vectors, $\beta_v$ and $\delta_v$, with element index $n$. Using the discrete Radon transform \eqref{eq:DiRadon} we define discretized versions of $B_\theta$ and $\phi_\theta$:
\begin{align}
B_v &= \phantom{-}\tfrac{2\pi}{\lambda} \mathbf{A}\beta_v \label{eqn:absorpv}, \\
\phi_v &= -\tfrac{2\pi}{\lambda} \mathbf{A}\delta_v \label{eqn:phaseshv}.
\end{align}
This leads to discrete versions of the transmittance and the Fresnel propagator, and from those the intensity $I^R_v$:
\begin{align}
T_v &= \exp(-B_v+\hat{\iota}\phi_v), \\
P^R_v &= -\tfrac{\hat{\iota}}{\lambda R}\exp\Big( \tfrac{\hat{\iota}\pi}{\lambda R}|t_v|^2 \Big), \\
I^R_v &= |T_v\star P^R_v|^2.
\end{align}
Here $\exp$ and $|\cdot|^2$ denote element-wise operations and $\star$ denotes the discrete convolution. In practice the discrete convolution is done in Fourier space.

\section{Reconstruction methods} \label{sec:RecMethods}

In this section we review the two methods that we compare, namely, a classical two-stage method and a combined method presented by Kostenko et al. Both methods are equipped with TV-regularization, which is described last.

\subsection{Two-stage method} \label{sec:TSmethod}
For absorption CT the measured intensity data can be directly related to the material properties, i.e., the attenuation coefficient, by Lambert-Beer's law \cite{Buzug2008}. For PCT the intensity measurements are related to the material properties through the non-linear propagation model Eqs.~(\ref{eq:Btheta})--(\ref{eqn:intensity}). The standard reconstruction process for PCT is a two-stage procedure, consisting of a \textit{phase retrieval} stage and a \textit{tomographic reconstruction} stage.

Different phase retrieval methods have been suggested and investigated in the literature. An introduction to, and a comparison of, some of them can be found in \cite{Burvall2011,Fienup1982,langer2008}. Following the work by Kostenko et al. \cite{Kostenko2013b,Kostenko2013a,Kostenko2013} we focus on the contrast transfer function (CTF) method, which is based on an assumption of low absorption, and which is derived from the expression of the intensity given in \eqref{eqn:intensity}. Based on a Taylor expansion of the transmittance function \eqref{eqn:transmit}, a CTF model that is linear in $B$ and $\phi$ can be derived in Fourier space \cite{Cloetens1999}.

For the case of a single detector distance and the so-called duality version of the CTF model with proportionality constant $\sigma = -\delta / \beta$, we have
\begin{align}
\widehat{I^R}(\omega) \approx\; &[ 2\sigma\sin(\pi\lambda R |\omega|^2)-2\cos(\pi\lambda R|\omega|^2)]\hat{B}(\omega) \nonumber \\ &+\delta_{\text{Dirac}}(\omega).  \label{eqn:CTFD}
\end{align}
Here $\delta_{\text{Dirac}}$ is the Dirac delta function and $\omega$ denotes the spatial frequency.
For $p$ being the physical size of a detector pixel, the sampling distance becomes $F_s = 1/p$ and hence $\omega \in [-\tfrac{F_s}{2},\tfrac{F_s}{2}]$. Introducing discrete frequency values $\omega_m$ the CTF method can be formulated as a discrete linear inverse problem: %
\begin{align} \label{eqn:CTFduality}
%\widehat{I^R_v} = \underbrace{\left[-2\mathbf{C}^R+2\sigma \mathbf{S}^R \right]}_{\mathbf{W}}\widehat{B}_v .
\widehat{I^R_v} &= \mathbf{W}_\text{TSD}\widehat{B}_v,\\
\mathbf{W}_\text{TSD} &= -2\mathbf{C}^R+2\sigma \mathbf{S}^R,
\end{align}
where $\mathbf{C}^R$ and $\mathbf{S}^R$ are diagonal matrices with $m$th diagonal elements $\cos{\psi_m}$ and $\sin{\psi_m}$, respectively, with
\begin{align}
\psi_m=\pi\lambda R|\omega_m|^2, \qquad m=1,2,\ldots,M.
\end{align}
%The resulting matrix system becomes
%
For the \textbf{two-stage} method based on the CTF \textbf{duality} method we use the name \textbf{TSD}.

The phase-retrieval stage of this method consists of solving \eqref{eqn:CTFduality} for $\hat{B}_v$, followed by an inverse Fourier transform and then solving \eqref{eqn:absorpv} for $\beta_v$. Since $\beta_v$ is assumed to be proportional to $\delta_v$ for this method, we can easily find $\delta_v$ afterwards.

Following the phase-retrieval stage is a stage where a CT reconstruction method is used to compute the object of interest. Typically CT reconstructions are carried out by using classical methods such as filtered backprojection (FBP) \cite{Buzug2008} or the algebraic reconstruction technique (ART) \cite{K1937}. 

\subsection{Algebraic combined method} \label{sec:ACmethod}

Methods that combine the phase-retrieval stage and the reconstruction stage have been proposed with different aims. In \cite{Bronnikov2002} a filtered backprojection type algorithm is derived and tested while in \cite{Kostenko2013a} an \textit{algebraic} combined model is presented and tested, showing promising preliminary results.
% are presented for this new approach and we therefore wish to make a thorough comparison between this algebraic combined method and the standard two-stage method.

In \cite{Kostenko2013a} Kostenko et al. also suggested that the standard phase-retrieval techniques could benefit from using the redundancy within an entire sinogram rather than just being based on the individual projections.

The term 'algebraic combined' in \cite{Kostenko2013a} refers to a combination of the linear operators, 
%This method was proposed in \cite{Kostenko2013a}.
$\mathbf{A}$, describing the discrete Radon transform, $\mathbf{F}$, the discrete Fourier transform, and $\mathbf{W}_\text{TSD}$ of the linear phase retrieval method, i.e.,
\begin{equation}
 \mathbf{W}_\text{ACD} = \tfrac{2\pi}{\lambda} \, \mathbf{W}_\text{TSD}\, \mathbf{F} \mathbf{A}.
\end{equation}
The discrete Fourier transform is introduced because
it is computationally convenient to
formulate the phase retrieval model as a matrix multiplication in Fourier space.
The algebraic combined method based on the duality version of the CTF method works on the resulting linear system
\begin{align} \label{eqn:ACD}
\widehat{I^R_v} = \mathbf{W}_\text{ACD} \beta_v \,.
\end{align}
This combined reconstruction method is called the \textbf{algebraic combined duality} method, \textbf{ACD}.

\subsection{Total variation regularization}

In absorption CT, TV-regularization has been shown to be advantageous for objects with piecewise constant material parameters \cite{Bian2010,Jorgensen2013}. TV-regularization preserves edges while smoothing away noise inside homogeneous regions. For the discrete linear problem \eqref{eq:DiRadon},
we formulate TV-regularization with regularization parameter $\alpha\in \mathbb{R}^+ $ as:
\begin{align} \label{eqn:TVform}
z_v^\alpha = \underset{u_v}{\text{argmin}} \left\{ \|\mathbf{A} u_v-b_v\|^2+\alpha\sum_{n=1}^{N^2} \|\mathbf{D}_n u_v \| \right\}.
\end{align}
Here $n$ is the pixel index, $N^2$ is the total number of pixels, assuming a square domain, and $\mathbf{D}_n u_v$ is the local finite difference gradient at pixel $n$.

Many objects from materials science which would be desirable to analyse with PCT have the property that they have approximately
piecewise constant material parameters, e.g., in the form of grains.
With this motivation Konstenko et al. proposed to incorporate TV-regularization into both the TSD and ACD methods, i.e., in \eqref{eqn:CTFduality} and \eqref{eqn:ACD}.

In the present work we also consider TV-regularization for both methods and in the remainder of the article by TSD and ACD we refer to the TV-regularized problems. In a direct assessment of the effect of combining linear operators in ACD, and not the effect of TV-regularization itself, we find it most appropriate to employ TV-regularization also in TSD. We note that one motivation for ACD is precisely the use of regularization which through the combination of linear operators regularizes the entire reconstruction problem including the phase-retrieval step. This is in contrast to the TV-regularized two-stage method where only the latter reconstruction step is regularized leaving the sensitive phase-retrieval stage unregularized. 

\section{Our contributions}

In \cite{Kostenko2013a} Kostenko et al. compared TV-regularized ACD and TSD and demonstrated improvements obtained by ACD in terms of root-mean-square error in most of their simulation experiments. We find that their pioneering results indicate a large potential for ACD, however we also point to several aspects in which the provided numerical evidence of reconstruction improvements by ACD may be improved:
\begin{enumerate}
 \item The positive results for ACD are for test images with one specific choice of material parameters that appears to not be motivated from physical materials. Thus, it remains open whether as clear improvements can be seen for test images with physical material parameters.
 \item The positive results for ACD are for noise-free data. In fact in a simulation study with noisy data, Kostenko et al. find the TSD to be superior. Only one noise level is considered so it remains unclear whether the combination of phase retrieval and reconstruction stages leads to a more noise-robust method.
 \item In their implementation of the optimization algorithm to solve the TV-regularized problem, Kostenko et al. \cite{Kostenko2013a} describe they stop the iterative algorithm when the relative change in the objective function value from one iteration to the next is smaller than $10^{-5}$. This is an intuitive choice, however well-known in the field of optimization to be heuristic and not guarantee closeness to the solution. This is because the iterative algorithm may occasionally take short steps while still far from the solution. This means that we cannot be sure that the shown reconstructions are indeed accurate TV-solutions but may be arbitrary intermediate images produced by the iterative algorithm. In fact, this problem might affect their conclusion that TSD is more robust to noise than ACD.
\end{enumerate}

In the present work, we address all of these three problems. Regarding the third problem, we implement in our optimization algorithm a stopping criterion that does ensure convergence to the TV-regularized reconstruction, thereby removing any doubt whether the numerical solution returned by the algorithm is in fact the sought-after TV-regularized solution.

We address the first and second problems by providing two sets of carefully designed simulation experiments. The first set compares TSD and ACD on test images with a range of physical material parameters, while the second compares TSD and ACD with respect to increasing amounts of noise.

To make the most direct and fair comparison between TSD and ACD, we employ TV-regularization for both and apply the same optimization algorithm with the same stopping criterion. 
%and we , the same simulated data and the same stopping criterion for the iterative algorithm.
%
Before proceeding to the results of the comparisons, we describe in the next section our implementation details.

\section{Implementation} \label{sec:implem}
The system matrix $\mathbf{A}$ in Eqs.\ (\ref{eqn:absorpv}) and (\ref{eqn:phaseshv}) is large and sparse. This means that it can often be stored in the memory of a standard modern laptop. For the ACD method the dense matrix $\mathbf{F}$ makes the combined system matrix dense, thus making it infeasible to store in memory. We circumvent the problem by a matrix-free implementation, in which the applications of the forward operator and its conjugate transpose are done without explicitly forming the matrices.

When solving the reconstruction problems Eqs.\ (\ref{eqn:absorpv}), (\ref{eqn:phaseshv}) and (\ref{eqn:ACD}) we impose TV-regularization \eqref{eqn:TVform}. Solving such large-scale problems requires efficient algorithms. We chose to implement the Chambolle-Pock (CP) algorithm \cite{Chambolle2011,Chambolle2011a} since it was shown to converge faster than, e.g., the FISTA method \cite{Chambolle2011a} when solving problems of the form \eqref{eq:DiRadon}. Moreover the CP algorithm can be well suited in the context of CT \cite{sidky2012}.
%When solving the reconstruction problems Eqs. (\ref{eqn:absorpv}), (\ref{eqn:phaseshv}) and (\ref{eqn:ACD}) we want to impose TV-regularization (as in \eqref{eqn:TVform}) on the solutions due to the discrete nature of many materials micro-structures. Solving large-scale TV-regularization problems requires efficient first-order algorithms. The Chambolle-Pock (CP) algorithm \cite{Chambolle2011,Chambolle2011a} was in \cite{sidky2012} shown to be well suited in the context of CT. Moreover the CP algorithm was shown to converge faster than, e.g., the FISTA method \cite{Chambolle2011a} when solving problems on the form \eqref{eq:DiRadon}. Therefore we implemented a matrix-free version to solve \eqref{eqn:TVform} for TSD and ACD.

Our implementation is mainly based on algorithm 4 in \cite{sidky2012}, modified by the adaptive parameter approach presented in algorithm 2 in \cite{Goldstein2013}. This modified approach introduces a primal residual $p^{(k)}$ and a dual residual $d^{(k)}$ for iteration $k$. As mentioned in \cite{Goldstein2013} these residuals can also be used to define a stopping criterion, since for the CP algorithm we have that

\begin{equation}
\lim_{k\to\infty} \|p^{(k)} \|^2 + \|d^{(k)} \|^2 = 0.
\end{equation}
We implemented a stopping criterion of the form
\begin{equation}
\label{eq:stoprule}
\|p^{(k)} \|^2 + \|d^{(k)} \|^2 < \tau \left( \|p^{(1)} \|^2 + \|d^{(1)} \|^2 \right)
\end{equation}
%modified this by normalizing with the primal-dual residual from the first iteration, such that we terminate the algorithm when
%\begin{equation}
%\frac{\|p^{(k)} \|^2 + \|d^{(k)} \|^2}{\|p^{(1)} \|^2 + \|d^{(1)} \|^2} < \tau,
%\end{equation}
for a user-defined tolerance $\tau$. In all of our numerical experiments $\tau$ was set to $10^{-6}$.

The reconstruction stage of TSD involves multiplication with $\mathbf{A}$ and
its transpose in each iteration.
The ACD method requires, in each iteration, additional FFTs and multiplications
with the diagonal matrix $\mathbf{W}$, but both operations are much less
computationally demanding than the multiplications with $\mathbf{A}$ and
hence the computational overhead in an ACD iteration, compared to TSD, is small.

%The ACD method will in each iteration be more computationally demanding than the TSD method, since the system matrix involves a Fourier transform. On the other hand, the fast Fourier transform (FFT) is very efficient, so in this matrix-free set-up the application of matrix $\mathbf{A}$ is the most dominant, with respect to computational time.

All implementations and simulations are carried out in \textsc{Matlab}, and the implemented code is available for download at \cite{PCTsimulator}. We also use the function \texttt{phantomgallery} from the \textsc{AIR Tools} package, Version 1.3 \cite{Hansen2012}, and the function \texttt{parbeam} from the Projector-Pack package, Version 0.2 \cite{ProjectorPack}.

%\todo{elaborate on Chambolle-Pock algorithm and implementation choices: adaptive parameter version, stopping criterion, etc.. \cite{Zabler2005},. Include a reference to appendix for testing the implementation.}

%%%%%%%%%%%%%%%%%%%%%%%%%%%%%%%%%%%%%%%%%%%%%%%%%%%

\section{Simulation results} \label{sec:comp_simu}

We compare the TSD and ACD methods across different material parameters and increasing noise levels. 
The comparisons rely on simulated data, carefully modelled to resemble data from a real physical set-up. Reconstructions are assessed in terms of achievable quality, compared to the ground truth and between the methods.

\subsection{Experimental set-up}
\begin{figure}[tbp]
\centering
\includegraphics[width=0.43\linewidth]{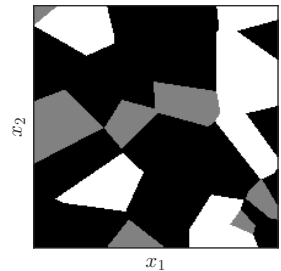}
\caption{2D phantom with a background material, black, and two different grain materials, gray and white.}
\label{fig:phantom}
\end{figure}

The simulated experiments are inspired by materials science where mappings of structures on a micrometer-scale are desired. For polycrystalline materials the structures are made up of \textit{grains}, and to mimic this we use a phantom which resembles grain-structure consisting of three different materials. The phantom, which is shown in \autoref{fig:phantom}, consists of one background material and grains of two different materials, all of which are described by indices $\beta$ and $\delta$. Indices of the used materials are listed in \autoref{tab:deltabeta}. If nothing else is mentioned, polycarbonate is used for the background material.

\renewcommand{\arraystretch}{1.2}

\begin{table}[tbp]
\centering
{\small
\caption{Absorption index $\beta$ and refractive index decrement $\delta$ for 40 keV X-rays for the simulated materials. From \cite{BetaDeltaCalc}. }
\begin{tabular}{l|l|l}
\hline 
\textbf{Material} & $\beta$ & $\delta$ \\
\hline
Polycarbonate (C$_{16}$H$_{14}$O$_3$) & $8.43\cdot 10^{-12}$ & $1.64\cdot 10^{-7}$ \\
\hline
Carbon (diamond) & $1.90\cdot 10^{-11}$ & $4.55\cdot 10^{-7}$ \\
\hline
Magnesium & $1.15\cdot 10^{-10}$ & $2.22\cdot 10^{-7}$ \\
\hline
Aluminium & $2.32\cdot 10^{-10}$ & $3.37\cdot 10^{-7}$ \\
\hline
Silicon & $2.68 \cdot 10^{-10}$ & $3.01\cdot 10^{-7}$ \\
\hline
Iron & $6.42\cdot 10^{-9}$ & $9.54\cdot 10^{-7}$ \\
\hline
Copper & $9.96\cdot 10^{-9}$ & $1.06\cdot 10^{-6}$ \\
\hline
\end{tabular}
  \label{tab:deltabeta} 
  }
\end{table}

% We compare the two reconstruction methods TSD and ACD, in terms of achievable quality of the reconstructed images. Both methods are used on the same simulated data, with the same optimization algorithm and compared using the same error measure and visualization in order to make as fair a comparison as possible.

Our set-up allows us to simulate free-space propagation PCT experiments and compare reconstruction methods. In our simulations we choose specific settings realistic for real physical experiments on a laboratory X-ray CT scanner. The chosen parameters are presented in \autoref{tab:exparam}.
% where $N_0$ represents the detected photons per pixel and per projection.
To make the simulated data more realistic, Poisson distributed noise is used to perturb the measured intensity.

In addition to the settings in \autoref{tab:exparam} the duality method requires a qualified guess on the proportionality constant $\sigma$ between $\beta$ and $\delta$. This has for all simulations been chosen as the exact proportionality for the grain material with the smallest $\beta$, e.g., for the first row in \autoref{fig:TSD_ACD_m} $\sigma = 	-1.95\cdot 10^{4}$. Experimental testing with different choices of $\sigma$ have shown that the impact of changing $\sigma$, within $\pm$20\% from the exact value, was negligible.

\begin{table}[tbp]
\centering
\caption{Parameters used in the simulations.}
\begin{tabular}{p{2.1cm} | p{5.5cm}}
\hline
\textbf{Parameters} & \textbf{Settings} \\
\hline
\textbf{Object:} & 2D and 200 x 200 pixels, pixel size 1 $\mu$m. \\
\textbf{X-ray} & Energy $40$keV. \\
\textbf{source:} & wavelength $\lambda = 0.31 $\r{A}. \\
\textbf{Photons:} & $N_0 = 10^5$ photons incident on object, average, per pixel, per projection. \\
\textbf{Distance:} & $R=0.5$m. \\
\textbf{Detector:} & Pixel size of $1\mu$m, $572$ pixels. \\
\textbf{Projections:} & 360 angles $\theta \in [0^\circ,180^\circ[\;$. \\
\hline
\end{tabular}
  \label{tab:exparam}
\end{table}

The regularization parameter was chosen empirically for each of the simulations in order to achieve the 'best' possible reconstruction. The 'best' reconstruction is in this work measured by two different means: A relative error measure
\begin{equation} \label{eqn:error}
E = \|u-u^*\| / \|u^*\|, \quad u^* = \text{original},
\end{equation}
and a visual comparison where sharp edges are favoured.
In the figures with the reconstructions we list the specific regularization parameter choices.

The reconstructions are visualized as images using grey-scale color-range $[0.9\cdot\min(u^*),1.1\cdot\max(u^*)]$; intensity values outside are truncated to this range.

\subsection{Effects of material properties}
\begin{figure}[tbp]
\begin{center}
\includegraphics[width=0.9\linewidth]{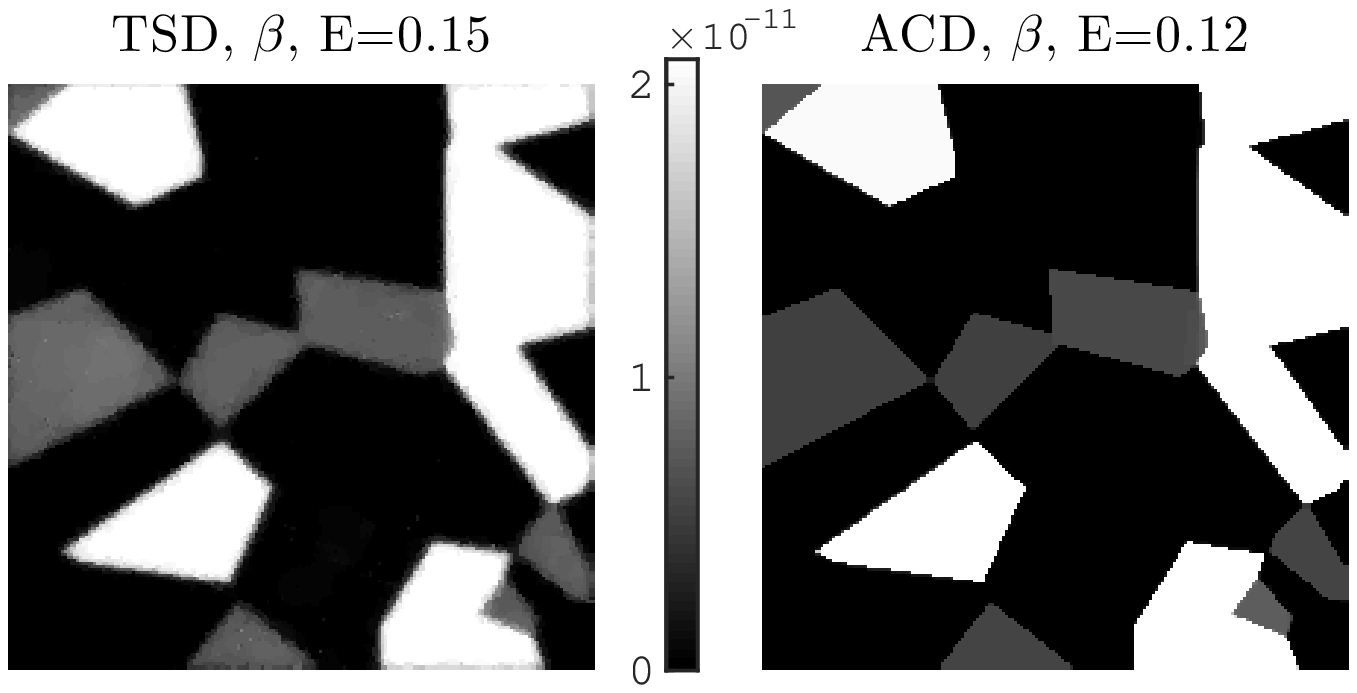}
\includegraphics[width=0.9\linewidth]{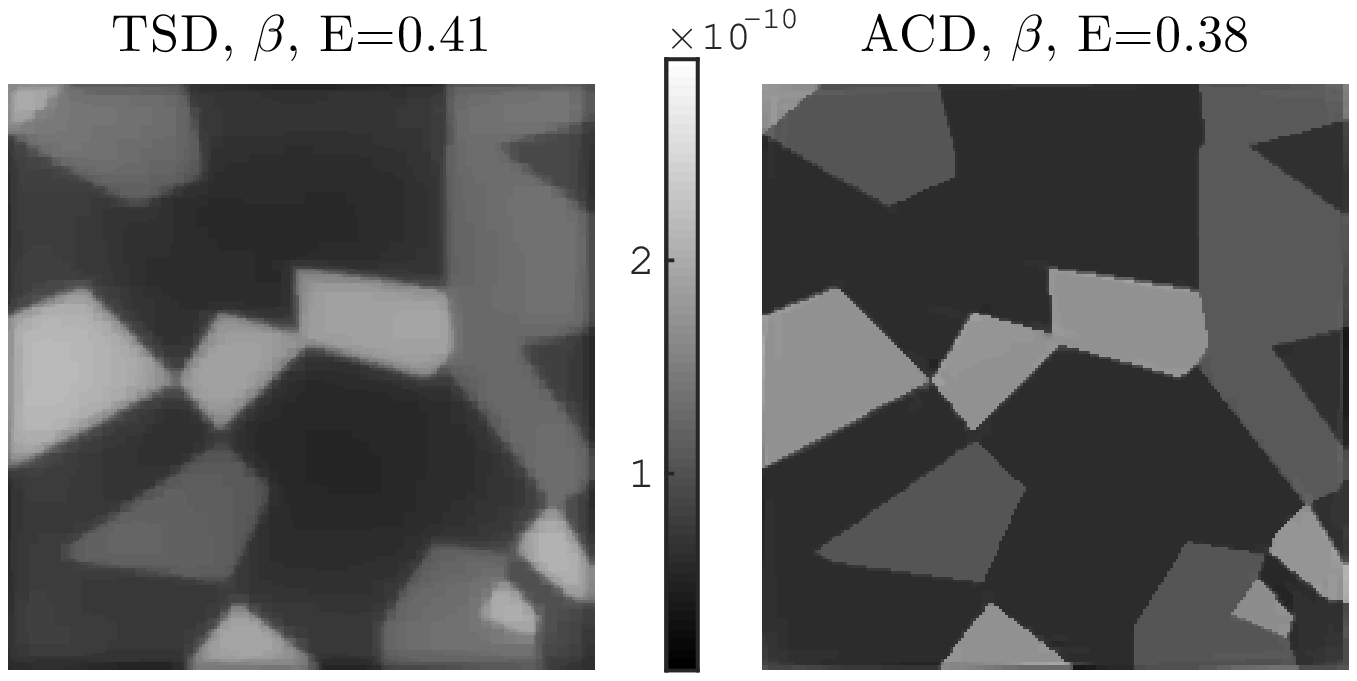}
\includegraphics[width=0.9\linewidth]{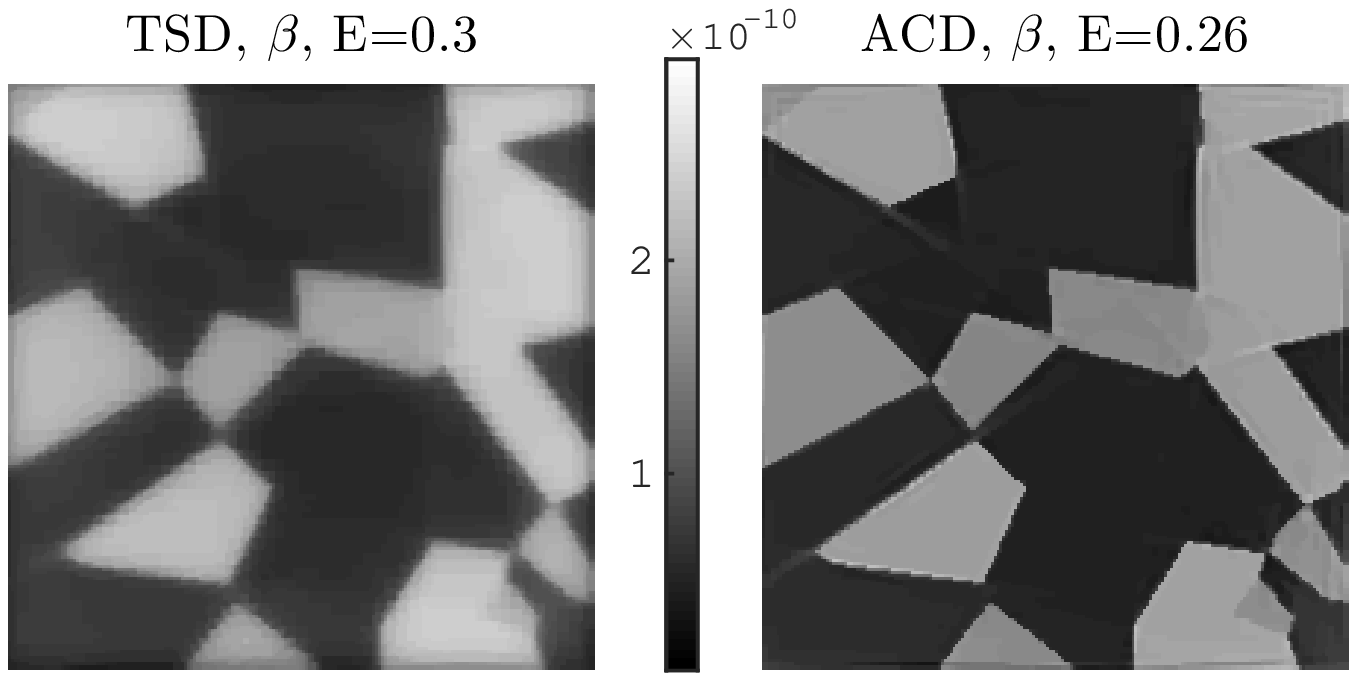}
\includegraphics[width=0.9\linewidth]{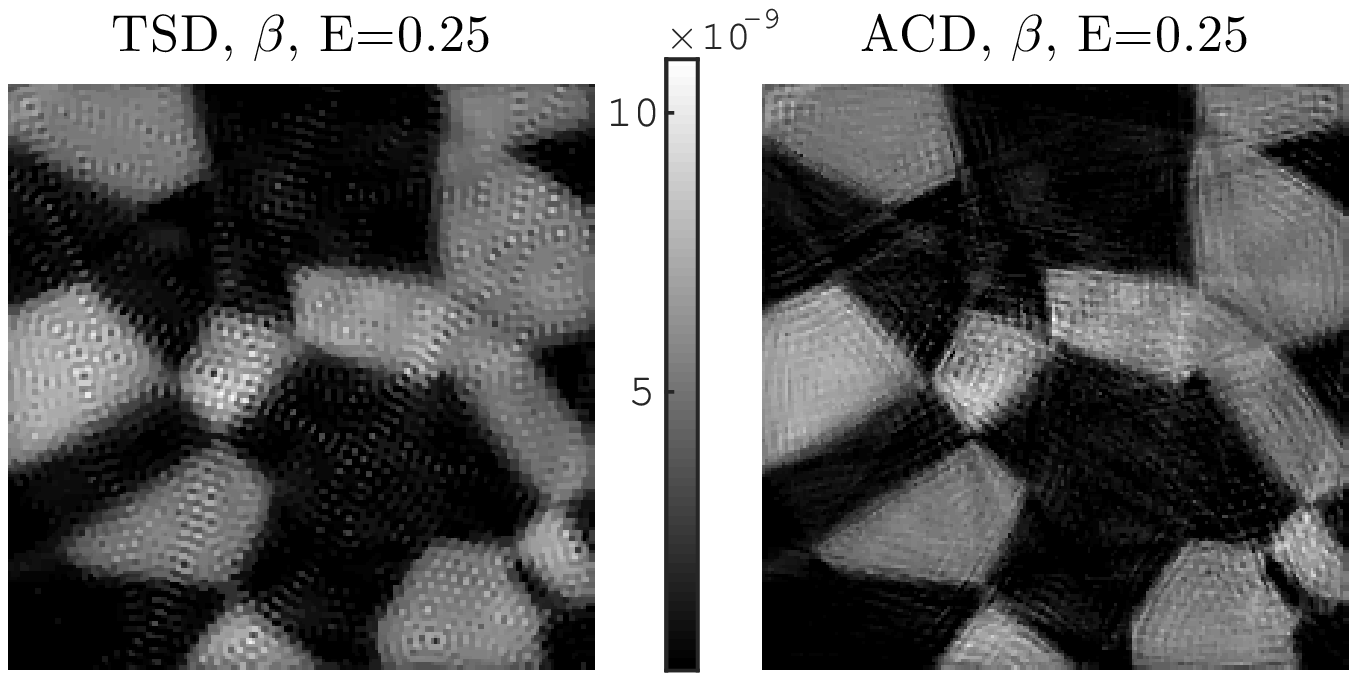}
\end{center}	
\caption{Simulations with different materials. First row: non-absorbing and non-refracting background, grains of polycarbonate and carbon in diamond form (highest $\beta$) and regularization parameters $\alpha_{TSD} = 0.01$, $\alpha_{ACD} = 120$. Second row: grains of silicon(highest $\beta$) and magnesium and regularization parameters $\alpha_{TSD} = 0.1$, $\alpha_{ACD} = 110$. Third row: grains of silicon (highest $\beta$) and aluminium and regularization parameters $\alpha_{TSD} = 0.12$, $\alpha_{ACD} = 112$. Fourth row: grains of copper (highest $\beta$) and iron and regularization parameters $\alpha_{TSD} = 6$, $\alpha_{ACD} = 5$. The error measure $E$ is defined in \eqref{eqn:error}.} \label{fig:TSD_ACD_m}
\end{figure}

The simulated phantom is varied with materials ranging from low-absorbing material to higher absorbing material, i.e., from low $\beta$ to higher $\beta$. The two-grain materials are chosen such that they have indices numerically close to each other since distinction between similar materials is the more challenging case in practical applications. Increasing the absorption will violate the low absorption assumption, which is part of the CTF model derivation, so higher absorption is also expected to increase the difficulty of the reconstruction problem. The reconstructions from our simulated experiments with materials of increasing absorption index are presented in \autoref{fig:TSD_ACD_m}.

For these four experiments the ACD method is generally seen to produce as good or better reconstructions than the TSD method, in terms of the error measure~$E$. The TSD results are all visually more blurry and with less sharp edges compared to the ACD results, even though both methods utilize the same TV-regularization method. In the process of choosing the 'best' reguralized reconstruction from a series of reconstructions (not shown here), it became clear that the TSD reconstructions were corrupted by artifacts and/or noise to a higher degree than the ACD reconstructions. On the four TSD reconstructions in \autoref{fig:TSD_ACD_m} this is what causes the reconstructions to be more blurred, since we gave more emphasis on the regularization term in order to compensate for noise and artifacts. We believe that the artifacts are due to the errors introduced by the linearization in the phase retrieval stage.

For the experiments with low absorption, in the first row of \autoref{fig:TSD_ACD_m}, distinction between the different materials is clear for both methods. The ACD reconstruction has sharper edges and a lower error measure than the TSD reconstruction.

For the silicon-magnesium and the silicon-aluminium experiments, in the second and third row, materials with similar chemical structures are seen to be harder to distinguish, as expected. The ACD method again produces reconstructions with sharper edges and a lower error-measure. For the silicon-aluminium reconstructions in row three, distinguishing between silicon and aluminium is difficult for both methods, and alternative methods using measurements from two more or more distances could improve these results.

In the experiments with high absorbing materials the reconstructions are highly affected by artifacts, such that distinction between background, grain, and artifact is difficult. In addition, the edges in the reconstructions are more blurry. Error measures does not tell the same story as the visual inspection, since they are relatively low compared to the
low-absorption experiments; this is due to the background material being correctly reconstructed.

The ACD method is computationally more demanding than the TSD method, because a larger number of iterations is needed to achieve the same solution accuracy. In the cases studied here, 1.4 -- 6.5 times more iterations were needed for the ACD method when using the stopping criterion in \eqref{eq:stoprule}.

\subsection{Effect of noise}
\begin{figure}[tbp]
\begin{center}
\includegraphics[width=0.9\linewidth]{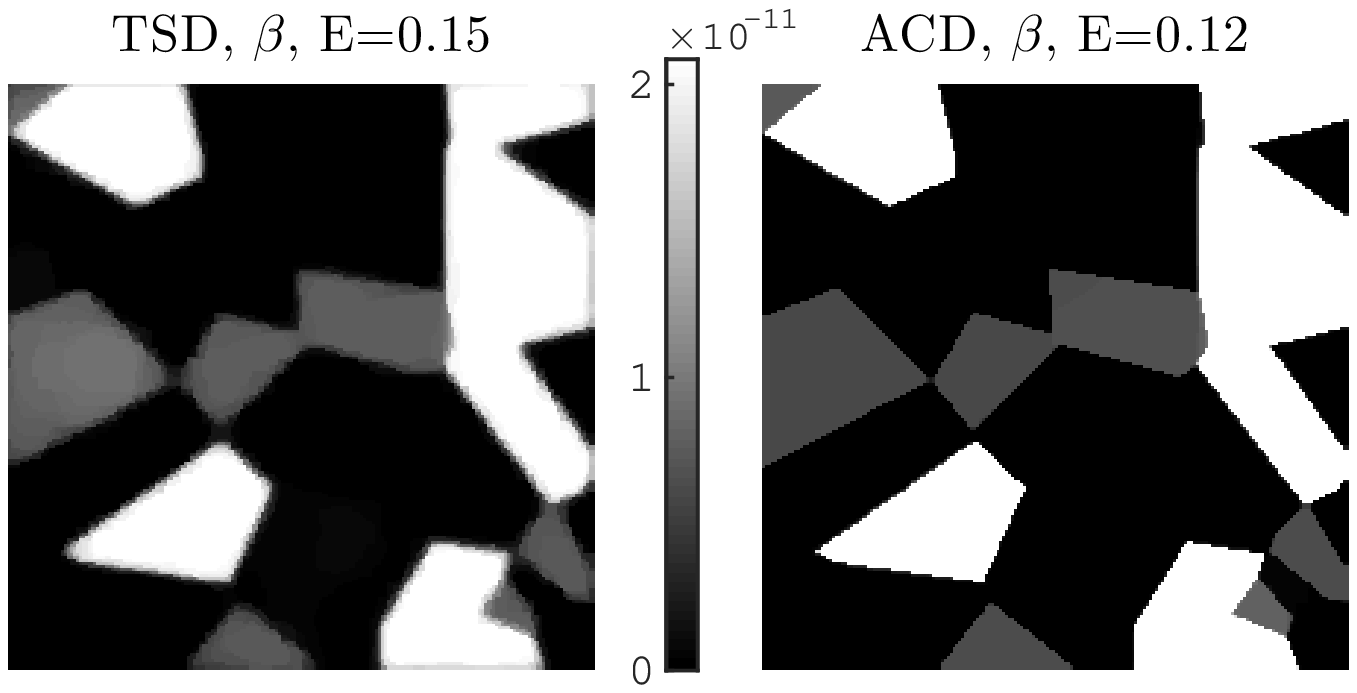}
\includegraphics[width=0.9\linewidth]{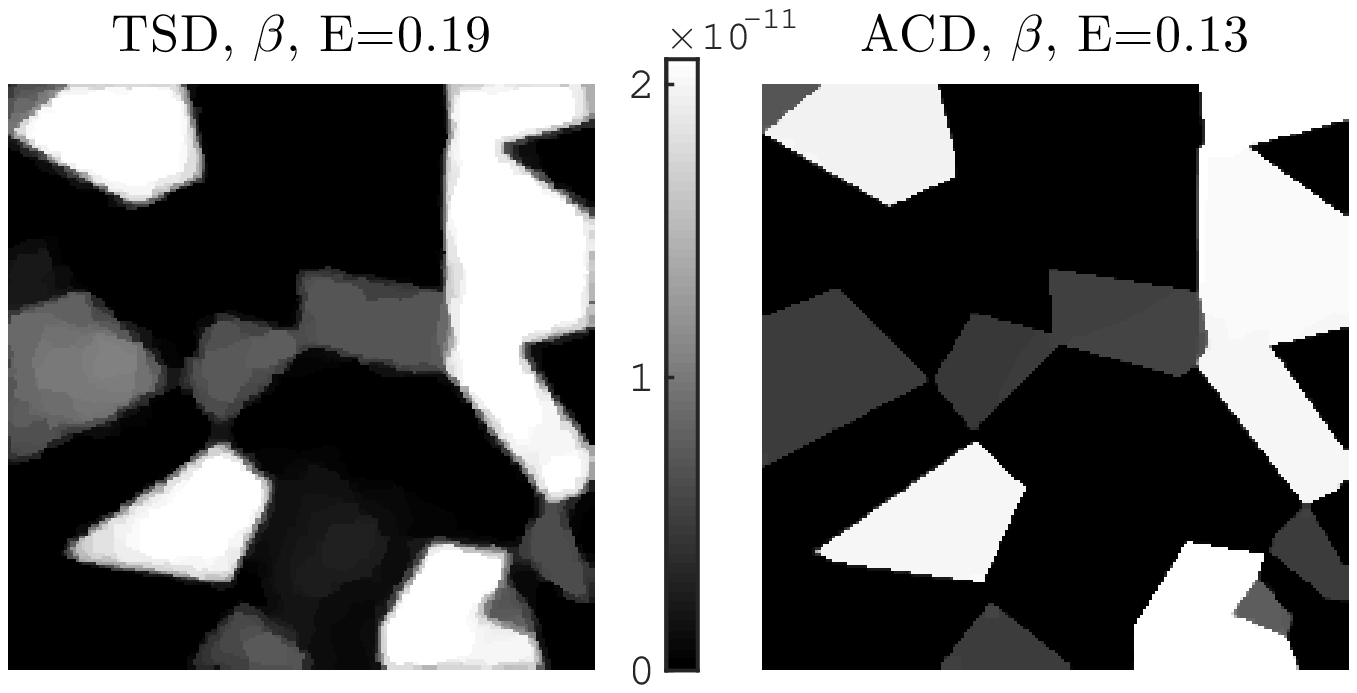}
\includegraphics[width=0.9\linewidth]{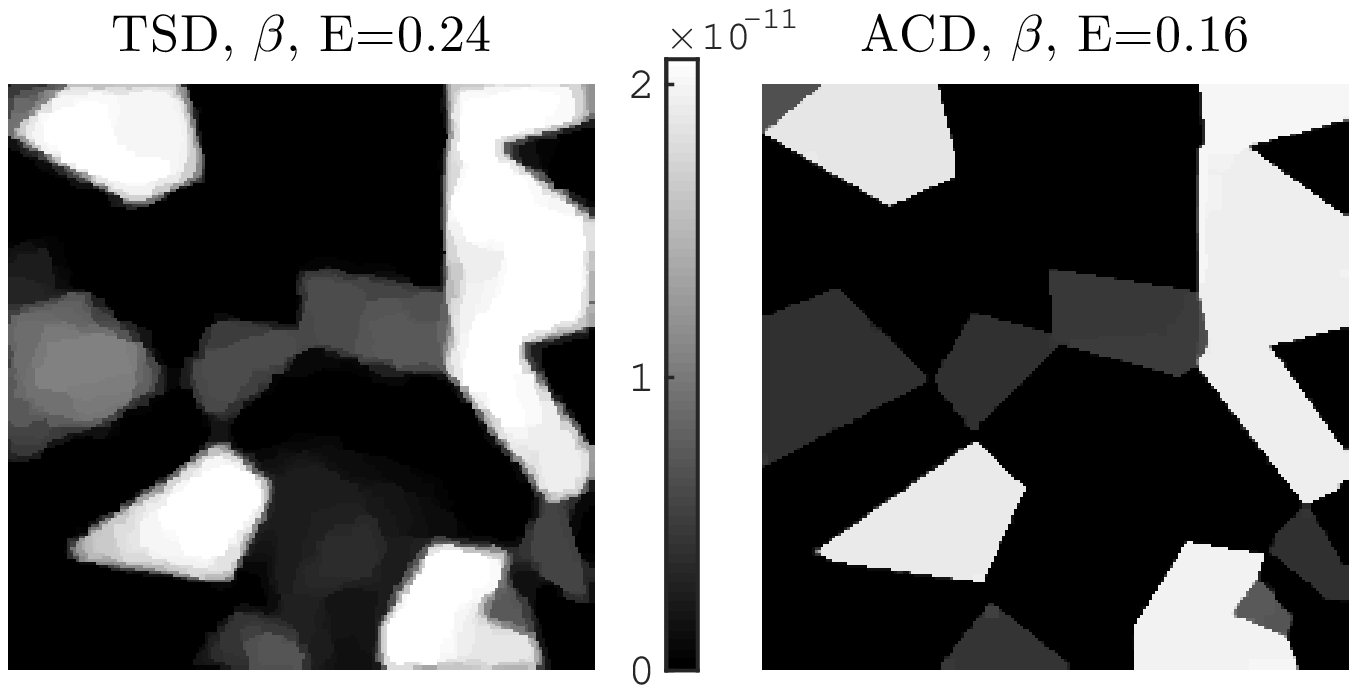}
\includegraphics[width=0.9\linewidth]{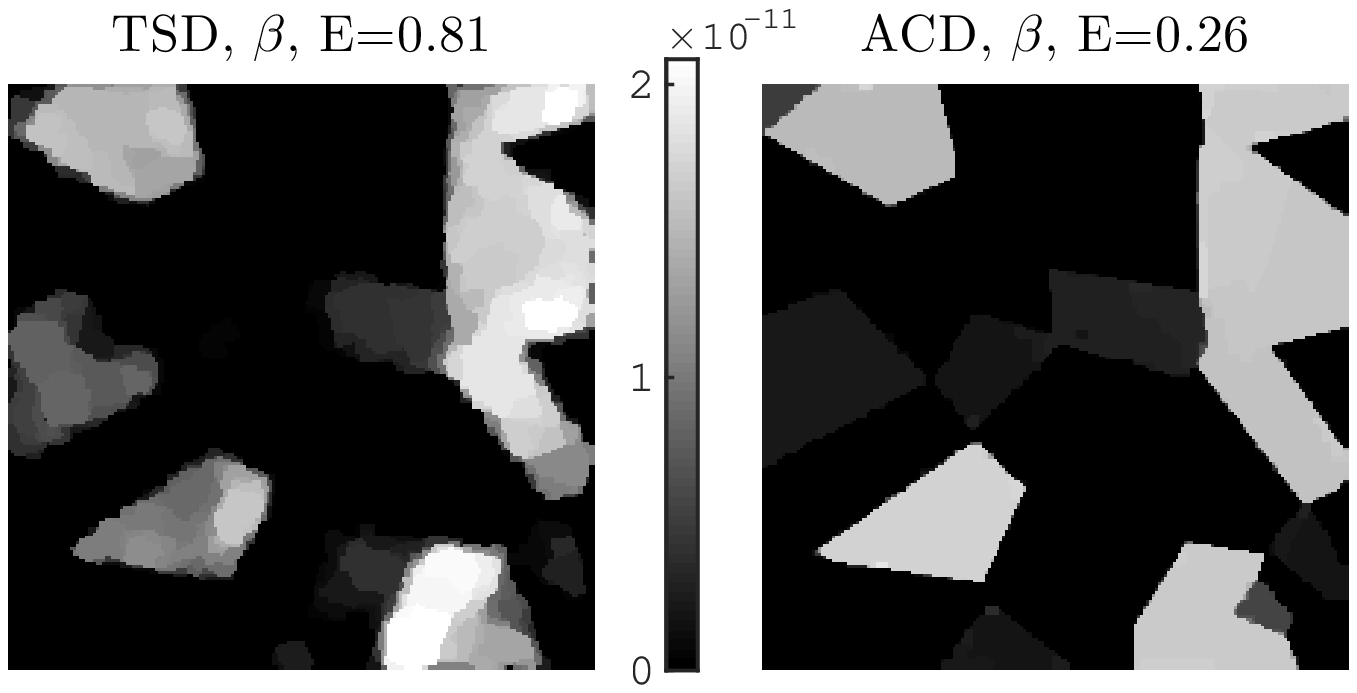}
\end{center}
\caption{Simulations for increasing data noise, i.e., decreasing number of photons $N_0$. Non-absorbing and non-refracting background. Grains of polycorbonate and carbon in diamond form (highest $\beta$). First row $N_0 = 5\times 10^{4} $, second row $N_0 = 10^{4} $, third row $N_0 = 5\times 10^{3} $ and fourth row $N_0 = 10^{3}$, photons per pixel. Regularization parameter $\alpha$ for the TSD method from top to botom were $[0.03,0.05,0.08,0.14]$ and for the ACD method $[100,140,170,300]$.} \label{fig:TSD_ACD_n}
\end{figure}

\begin{figure}[tbp]
\centering
\includegraphics[width=0.8\linewidth]{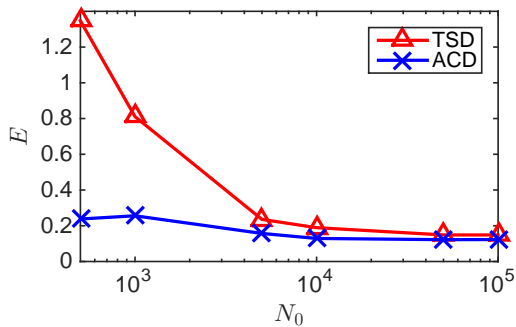}
\caption{Relative reconstruction error $E$ for increasing number of photons $N_0$, i.e., decreasing noise in the data.}
\label{fig:c_TSD_ACD_e}
\end{figure}

\begin{figure}[tbp]
\centering
\includegraphics[width=0.8\linewidth]{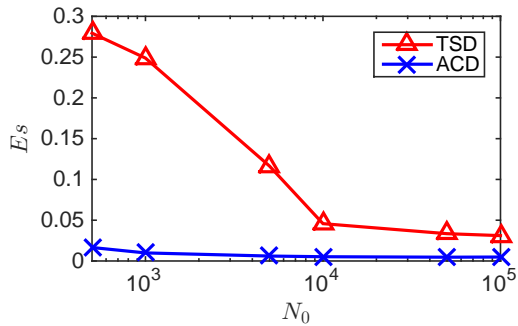}
\caption{Segmentation error $Es$ for increasing number of photons $N_0$, i.e., decreasing noise in the data.}
\label{fig:c_TSD_ACD_es}
\end{figure}

The reconstructions from simulated experiments with a decreasing number of recorded photons, and hence increasing noise levels, are presented in \autoref{fig:TSD_ACD_n}. The error measure $E$ is plotted for increasing $N_0$ (cf.\ Table \ref{tab:exparam}) in \autoref{fig:c_TSD_ACD_e}. Using Otsu's simple thresholding segmentation method \cite{Otsu1979} on the reconstructions, the segmentation errors
\begin{equation}
Es = \#\text{misclassified pixels} / \#\text{pixels}
\end{equation}
are calculated and plotted against $N_0$ in \autoref{fig:c_TSD_ACD_es}.

The TSD reconstructions are seen to deteriorate as $N_0$ decreases (and the relative noise increases), where the grains of the lowest absorbing material closest to the object center become indistinguishable from the background -- cf.\ the bottom row in \autoref{fig:TSD_ACD_n}. The edges become more blurry and misclassification of the grains is likely to occur. The error measure increases drastically to a limit where the reconstructions are unreliable.

ACD reconstructions show much greater robustness to the noise: edges remain sharp, materials can be distinguished, and the error varies slowly with $N_0$. For the problems with higher relative noise (smaller $N_0$), the polycarbonate grains can be visually hard to distinguish from the background in the chosen gray-scale -- though numerically the difference is still distinct as validated by the low segmentation errors in \autoref{fig:c_TSD_ACD_es}.

% \figoneF{c_TSD_ACD_e}{0.45}{Relative reconstruction error $E$ for increasing number of photons $N_0$, i.e., decreasing noise in the data.}

%%%%%%%%%%%%%%%%%%%%%%%%%%%%%%%%%%%%%%%%%%%%%%%%%%%

\FloatBarrier

\section{Conclusion}

The simplicity of the geometry in the free-space propagation, one sample-detector distance approach makes it attractive for many experiments, including those where speed of data acquisition or dose is a limitation.  Likewise, experience has shown that total variation (TV) regularization works well for absorption or phase contrast tomography on a large class of materials comprising disjunct phases, cracks, pores, etc.

The outcome of the simulations performed here is that the number of photons required to compute a reconstruction of a certain quality can be reduced substantially and the combined reconstruction method is therefore of general interest to the X-ray (and neutron) imaging community. We emphasize that the suggested combined method is more computationally demanding than the classical two-stage method, and hence less suitable for on-line or real-time processing.
\section*{Funding Information}

This work is part of the project HD-Tomo funded by Advanced Grant No.\ 291405 from the European Research Council. Henning Friis Poulsen acknowledges the European Research Council Advanced Grant "Diffraction-based Transmission Microscopy" No. 291321.\\

\section*{Acknowledgments}

The authors thank A. Kostenko and J. Batenburg insightful discussions and for sharing their work and code on phase-contrast tomography. We also thank M. S. Andersen for discussions related to numerical optimization and Y. Dong for discussions about image processing in general.

% Bibliography
\bibliographystyle{abbrv}
\bibliography{MasterThesis,links}

\end{document}